\documentclass[
  ,draft            % use draft while you are working on the paper
  ]
  {aipproc}
\usepackage{amssymb}
\layoutstyle{8x11single}

\newtheorem{theorem}{Theorem}

\newcommand{\qed}{\nobreak \ifvmode \relax \else
      \ifdim\lastskip<1.5em \hskip-\lastskip
      \hskip1.5em plus0em minus0.5em \fi \nobreak
      \vrule height0.75em width0.5em depth0.25em\fi}

\newcommand{\re}{\mathrm{Re}}% partie reelle
\newcommand{\im}{\mathrm{Im}}% partie imaginaire
\newcommand{\Log}{\mathrm{Log}}% logarithme principal
\newcommand{\Li}{\mathrm{Li}}% polylogarithme

\begin{document}

\title{On Two Applications of Herschel's Theorem}



\author{Lazhar Fekih-Ahmed}{
  address={\'{E}cole Nationale d'Ing\'{e}nieurs de Tunis, BP 37, Le Belv\'{e}d\`{e}re 1002 , Tunis, Tunisia}
}

\begin{abstract}
As a first application of a very old theorem, known as Herschel's
theorem, we provide direct elementary proofs of several explicit
expressions for some  numbers and polynomials that are known in
combinatorics. The second application deals with the analytical
continuation of the polylogarithmic function of complex argument
beyond the circle of convergence.
\end{abstract}

\maketitle

%%%%%%%%%%%%%%%%%%%%%%%%%%%%%%%%%%%%%%%%%%%%
%% MAINMATTER
%%%%%%%%%%%%%%%%%%%%%%%%%%%%%%%%%%%%%%%%%%%%

\section{Herschel's Theorem}\label{sec1}

Herschel's theorem gives the expression of the expansion of  a
function of the form $\phi(e^{-t})$ into a Taylor series in
ascending powers of $t$. The proof of Herschel's theorem is
straightforward: By Taylor's theorem, we have

\begin{eqnarray}
\phi(e^{-t})&=&\phi(1-(1-e^{-t}))=\phi(1)+\phi^{\prime}(1)\frac{(-1)}{1!}(1-e^{-t})+
\phi^{\prime\prime}(1)\frac{(-1)^2}{2!}(1-e^{-t})^2+\nonumber\\
&&\cdots+\phi^{(n)}(1)\frac{(-1)^{n}}{n!}(1-e^{-t})^n+\cdots\nonumber\\
&=&\sum_{n=0}^{\infty}\phi^{(n)}(1)\frac{(-1)^{n}}{n!}(1-e^{-t})^n.\label{sec1-eq1}
\end{eqnarray}

If we let

\begin{equation}\label{sec1-eq2}
\phi(e^{-t})=a_{0}+a_{1}t+a_{2}t^2+\cdots+ a_{n}t^n+\cdots
=\sum_{n=0}^{\infty}a_{n}t^n,
\end{equation}

then the coefficient $a_n$ will be equal to the sum of the
coefficients of $t^n$ in the expansion of all the terms in the
right hand side of equation~(\ref{sec1-eq1}). Now, we know that

\begin{equation}\label{sec1-eq3}
(-1)^{k}(1-e^{-t})^k=(-1)^{k}{k\choose
0}e^{-0.t}+(-1)^{k-1}{k\choose 1}e^{-t}+\cdots-{k\choose
k-1}e^{-(k-1)t}+{k\choose k}e^{-kt},
\end{equation}

and so the coefficient of $t^n$ in the right hand side of
equation~(\ref{sec1-eq3}) is equal to

\begin{equation}\label{sec1-eq4}
\frac{(-1)^{n}}{n!}\bigg[(-1)^{k}{k\choose
0}0^{n}+(-1)^{k-1}{k\choose 1}1^{n}+\cdots -{k\choose
k-1}(k-1)^{n}+{k\choose k}k^{n}\bigg].
\end{equation}

But using the notation of the calculus of finite differences, the
last equation can be written as

\begin{equation}\label{sec1-eq5}
\frac{(-1)^{n}}{n!}\Delta^{k}0^{n}.
\end{equation}

Combination of the above equations leads to a formula for the
coefficient $a_n$, better known in the old literature as
Herschel's theorem \citep[chap.~2]{boole:2009},
\citep{Hamilton:1837}:

\begin{theorem}[Herschel's Theorem]\label{sec1-thm1}
\begin{equation}\label{sec1-eq6}
a_{n}=\frac{(-1)^{n}}{n!}\bigg[\phi(1).0^{n}+\frac{\phi^{\prime}(1)}{1!}\Delta
0^{n}+\frac{\phi^{\prime\prime}(1)}{2!}\Delta^{2}0^{n}+\cdots+
\frac{\phi^{(n)}(1)}{n!}\Delta^{n}0^{n}\bigg]=\frac{(-1)^{n}}{n!}
\sum_{j=0}^{n}\frac{\phi^{(j)}(1)}{j!}\Delta^{j}0^{n}.
\end{equation}
\end{theorem}

\section{Finding Explicit formulas using Herschel's Theorem}

Several explicit formulas can be deduced from Herschel's Theorem.
Let's suppose that we have a series of numbers  defined by a
generating function for which we want to find an explicit formula.
The first step consists in expressing the function as a function
of $e^{-t}$. The second step is to use Theorem~\ref{sec1-thm1} to
provide the explicit formula using finite differences of $0$. We
finish the sections  with some examples.

\subsection{Bernoulli Numbers}

The generation function of Bernoulli numbers is given by \cite[p.
48]{Comtet:1974}:

\begin{equation}\label{sec2-eq1}
\phi(e^{-t})=\frac{t}{e^{t}-1}=\frac{-\ln(e^{-t})e^{-t}}{1-e^{-t}}=\sum_{n=0}^{\infty}B_n\frac{t^n}{n!};
|t|<2\pi
\end{equation}

By setting $X=1-e^{-t}$, we have

\begin{eqnarray}\nonumber
\phi(1-X)&=&\frac{\ln(1-X)}{X}(1-X)=\frac{1}{X}\bigg[X+\frac{X^2}{2}+\frac{X^3}{3}+\cdots\bigg](1-X),\\
&=&1-\frac{X}{1.2}-\frac{X^2}{2.3}-\cdots-\frac{X^n}{n.(n+1)}-\cdots\label{sec2-eq2}
\end{eqnarray}

where in this case $\frac{\phi^{(j)}(1)}{j!}= -\frac{1}{j(j+1)}$
for $j\ge 1$. Herschel's theorem gives the following well-known
explicit formula for Bernoulli numbers:

\begin{equation}\label{sec2-eq3}
B_n=(-1)^{n}n!a_n=1.0^{n}+(-1)^{n+1}\frac{\Delta
0^{n}}{1.2}+(-1)^{n+1}\frac{\Delta^{2}0^{n}}{2.3}+\cdots
+(-1)^{n+1}\frac{\Delta^{n}0^{n}}{n(n+1)}.
\end{equation}

\subsection{Euler Polynomials}

Euler polynomials of degree $n$ in $x$ are denoted by $E_n(x)$ and
are defined by the generating function

\begin{equation}\label{sec2-eq4}
\phi(e^{-t})=\frac{2e^{tx}}{e^{t}+1}=\frac{2e^{-t(1-x)}}{1+e^{-t}}=\sum_{n=0}^{\infty}E_n(x)\frac{t^n}{n!};
|t|<\pi
\end{equation}

Again by setting $X=1-e^{-t}$, we rewrite the generating function
as

\begin{equation}\label{sec2-eq5}
\phi(1-X)=\frac{(1-X)^{1-x}}{1-\frac{X}{2}}.
\end{equation}

Since $|X|<1$, the generalized binomial theorem gives

\begin{eqnarray}\label{sec2-eq6}
(1-X)^{1-x}&=&\sum_{k=0}^{\infty}(-1)^k{1-x\choose
k}X^k=1+(x-1)X+\frac{(x-1)x}{2!}X^2+\frac{(x-1)x(x+1)}{3!}+\cdots\\
\frac{1}{1-\frac{X}{2}}&=&1+X+\frac{X^2}{2^2}+\frac{X^3}{3^2}+\cdots,\label{sec2-eq7}
\end{eqnarray}

The product of the two series provides

\begin{equation}\label{sec2-eq8}
\phi(1-X)=\sum_{n=0}^{\infty}\frac{1}{2^n}\sum_{k=0}^{n}\frac{(-1)^k}{2^k}{1-x\choose
k}X^n.
\end{equation}

In this case, $\frac{\phi^{(j)}(1)}{j!}=
\frac{(-1)^j}{2^j}\sum_{k=0}^{j}\frac{(-1)^k}{2^k}{1-x\choose k}$,
and therefore, Herschel's theorem gives the desired explicit
formula for Euler polynomials:

\begin{equation}\label{sec2-eq9}
E_n(x)=(-1)^{n}n!a_n=(-1)^{n}\sum_{j=0}^{n}\frac{(-1)^j}{2^j}\sum_{k=0}^{j}\frac{(-1)^k}{2^k}{1-x\choose
k}\Delta^{j}0^{n}.
\end{equation}

The reader can compare with the formula  obtained in
\cite{Luo:2006}. Clearly, one can also obtain an explicit formula
for Euler polynomials of higher order.

\subsection{Eulerian Numbers and Polynomials}

It is known that the classical Eulerian polynomials
$A_n(\lambda)$, $0<\lambda<1$ have the exponential generating
function \citep[p.51]{Comtet:1974}, \citep{Butzer-Hauss:1993}

\begin{equation}\label{sec2-eq10}
1+\sum_{n=1}^{\infty}\frac{A_n(\lambda)}{\lambda}\frac{t^n}{n!}=\frac{1-\lambda}{
e^{t(\lambda-1)}-\lambda},
\end{equation}

By replacing $t$ by $\frac{t}{\lambda-1}$ in (\ref{sec2-eq10}), we
obtain the following function, \citep[chap.~6,
p.244]{Comtet:1974}:

\begin{equation}\label{sec2-eq11}
1+\sum_{n=1}^{\infty}\frac{A_n(\lambda)}{\lambda(\lambda-1)^n}\frac{t^n}{n!}=\frac{1-\lambda}{e^{t}-\lambda},
\end{equation}

Carlitz~\citep{Carlitz:1959} denoted the numbers
$\frac{A_n(\lambda)}{\lambda(\lambda-1)^n}$ by $H_n(\lambda)$.

A theorem of Frobenius~\citep{Frobenius:1910} states that the
Eulerian polynomials are given by

\begin{eqnarray}
 A_n(\lambda) &=& \lambda\sum_{j=1}^{n}j!S(n,j)(\lambda-1)^{n-j} \label{sec2-eq12}\\
  &=& \lambda\sum_{j=1}^{n}(\lambda-1)^{n-j}\Delta^{j}0^{n},\label{sec2-eq13}
\end{eqnarray}
where $S(n,j)$ are the Stirling numbers of the second kind. Note
that Stirling numbers of the second kind are defined by

\begin{equation}\label{sec2-eq14}
S(n,j)=\frac{1}{j!}\sum_{k=0}^{j}(-1)^{j-k}{j\choose k}k^n,
\end{equation}
which can be easily written as a function of $k$'th forward
difference of $0^n$:

\begin{equation}\label{sec2-eq15}
S(n,j)=\frac{1}{j!}\Delta^{j}0^{n}.
\end{equation}

We now give a new formula for the Eulerian polynomials using
Herschel's theorem. The formula complements  Frobenius formula. We
start by setting $X=1-e^{-t}$ and rewriting the generating
function (\ref{sec2-eq11}) as

\begin{equation}\label{sec2-eq16}
\phi(e^{-t})=\frac{1-\lambda}{e^{t}-\lambda}=\frac{(1-\lambda)e^{-t}}{1-\lambda
e^{-t}}=\sum_{n=0}^{\infty}H_n(\lambda)\frac{t^n}{n!};
|t|<\log{\lambda}+2\pi.
\end{equation}

In terms of the variable $X$, the generating function becomes

\begin{equation}\label{sec2-eq17}
\phi(1-X)=\frac{1-X}{1-\frac{\lambda}{\lambda-1}X}=1+\sum_{n=1}^{\infty}
\frac{\lambda^{n-1}}{(\lambda-1)^{n}}X^n.
\end{equation}

By Herschel's theorem, we obtain $H_0(\lambda)=1$ and for $n\ge 1$

\begin{equation}\label{sec2-eq18}
H_n(\lambda)=(-1)^{n}\sum_{j=1}^{n}(-1)^j\frac{\lambda^{j-1}}{(\lambda-1)^{j}}\Delta^{j}0^{n}.
\end{equation}

Thus, $A_0(\lambda)=1$ and for $n\ge 1$
\begin{equation}\label{sec2-eq19}
A_n(\lambda)=\sum_{j=1}^{n}\lambda^{j}(1-\lambda)^{n-j}\Delta^{j}0^{n}.
\end{equation}

Note that equation (\ref{sec2-eq19}) reminds of Bernstein
polynomials with  coefficients as functions of differences of
zero.

\subsection{Genocchi Numbers}

Genocchi numbers are defined by the generating function \cite[p.
49]{Comtet:1974}.

\begin{equation}\label{sec2-eq20}
\frac{2t}{e^{t}+1}=\sum_{n=1}^{\infty}G_n\frac{t^n}{n!}.
\end{equation}

Using the change of variable $X=1-e^{-t}$, the generating function
becomes

\begin{equation}\label{sec2-eq21}
\phi(e^{-t})=\phi(1-X)=-\frac{(1-X)\log(1-X)}{1-\frac{X}{2}}
\end{equation}

\begin{eqnarray}\label{sec2-eq22}
\log(1-X)(1-X)&=&X-\frac{X^2}{1.2}-\frac{X^3}{2.3}-\cdots-\frac{X^{n}}{(n-1).n}-\cdots,\\
\frac{1}{1-\frac{X}{2}}&=&1+\frac{X}{2}+\frac{X^2}{2^2}+\cdots+\frac{X^n}{2^n}-\cdots\label{sec2-eq23}
\end{eqnarray}

Multiplying the two power series, the coefficient of $X^0$ is
zero, the coefficient of $X$ is equal to 1 and the coefficient of
$X^n$ is given by

\begin{equation}\label{sec2-eq24}
(-1)^n\frac{\phi^{(n)}(1)}{n!}=\frac{1}{2^{n-1}}\bigg[1-\sum_{k=2}^{n}\frac{2^{k-1}}{(k-1)k}
\bigg], n\ge 2.
\end{equation}

Finally, an application of Herschel's theorem yields the formula
for the Genocchi numbers:

\begin{equation}\label{sec2-eq25}
G_n=(-1)^{n}n!a_n=(-1)^{n}\bigg\{ 1+
\sum_{j=2}^{n}\frac{1}{2^{j-1}}\bigg[
1-\sum_{k=2}^{n}\frac{2^{k-1}}{(k-1)k}\bigg]\Delta^{j}0^{n}
\bigg\}.
\end{equation}

\section{An Analytic continuation of the Polylogarithm}

The polylogarithm $\Li_{s}(x)$ is defined by the power series

\begin{equation}\label{sec3-eq1}
\Li_{s}(x)=\sum_{n=1}^{\infty}\frac{x^n}{n^{s}}.
\end{equation}

The above definition is valid for all complex values $s$ and all
complex values of $x$ such that $|x|<1$.

The conformal mapping

\begin{eqnarray}\label{sec3-eq2}
x&=&1-e^{-t},\\
t&=&-\Log(1-x)\label{sec3-eq3}
\end{eqnarray}

maps the part of the $x$-plane between two half-lines starting
from the point $(1,0)$ to a strip parallel to the $x$-axis in the
$t$ plane. Moreover, the function $1-e^{-t}$ is conformal at each
point $t\in \mathbb{C}$, since its derivative does not vanish at
$t$. Its restriction to the horizontal strip $\{|\im(t)|<\pi \}$
is a conformal mapping of the strip onto the cut plane
$\mathbb{C}\setminus [1, \infty)$.

The principal branch $\Log(1-x)$ of $\log(1-x)$ is also a
conformal mapping of the cut plane $\mathbb{C}\setminus [1,
\infty)$ onto the horizontal strip $\{|\im(t)|<\pi \}$.

It is easy to verify that $\Li_{s}(x)$ has one finite singularity,
namely the point $x=1$. The point $x=1$ is mapped to $\infty$ by
the conformal mapping (\ref{sec3-eq3}). Making the substitution
(\ref{sec3-eq2}) into (\ref{sec3-eq1}), the series becomes

\begin{equation}\label{sec3-eq4}
\Li_{s}(1-e^{-t})=\sum_{n=1}^{\infty}\frac{(1-e^{-t})^n}{n^{s}}.
\end{equation}

Expanding the right hand side of (\ref{sec3-eq4}) into a series in
powers of $t$, we get

\begin{equation}\label{sec3-eq5}
\Li_{s}(1-e^{-t})=\sum_{n=0}^{\infty}a_n t^n.
\end{equation}

where the $a_n$ may be calculated using Herschel's theorem.
Indeed, we have here

\begin{equation}\label{sec3-eq6}
(-1)^n\frac{\phi^{(n)}(1)}{n!}=\frac{1}{n^{s}}, n\ge 1; \phi(1)=0.
\end{equation}

By Herschel's theorem, we thus have

\begin{equation}\label{sec3-eq7}
a_0=0,\quad a_n=\frac{(-1)^{n}}{n!}
\sum_{j=1}^{n}\frac{(-1)^{j}}{j^s}\Delta^{j}0^{n},\quad n\ge 1.
\end{equation}

Now if the function $\Li_{s}(x)$ is regular in the plane
$\mathbb{C}$ minus the  semiaxis $\re(x)>1$, the series
(\ref{sec3-eq5}) is necessarily convergent in the circle $|t|<1$.
Conversely, if the series (\ref{sec3-eq5}) is convergent  in the
circle $|t|<1$, then the function $\Li_{s}(x)$ is regular in the
cut plane $\mathbb{C}\setminus [1, \infty)$. Therefore, we can
assert that $\Li_{s}(x)$ can be represented at any point of the
cut plane $\mathbb{C}$ by the following expansion

\begin{equation}\label{sec3-eq8}
\Li_{s}(x)=\sum_{n=1}^{\infty}\bigg(  \frac{(-1)^{n}}{n!}
\sum_{j=1}^{n}\frac{(-1)^{j}}{j^s}\Delta^{j}0^{n}\bigg)(\Log(1-x))^n.
\end{equation}

There exists other integral and series relations which provide the
analytic continuation of the polylogarithm beyond the circle of
convergence $|x| = 1$ of the defining power series. But these
relations are valid for all but some exceptional values of $s$,
\citep[p.~139-140]{lindelof:1905},
\citep{jonquiere:1889,wirtinger:1905,Truesdell-Bateman:1945}. To
the author's knowledge, the series (\ref{sec3-eq8}) is the only
series that defines the polylogarithm for all values of $s\in
\mathbb{C}$ and all values of $x\in\mathbb{C}$.

\bibliographystyle{aipproc}

\end{document}